\documentclass[accepted,3p,10pt]{elsarticle}
\usepackage{amsmath}
\usepackage[mathlines]{lineno}
\usepackage{amssymb}
\usepackage{amsmath}
\usepackage{cases}
\usepackage{epstopdf}
\usepackage[caption=false,font=footnotesize]{subfig}
\usepackage{fixltx2e}
\usepackage{amssymb}
\usepackage{mathtools}
\usepackage{bm}
\usepackage{multirow}
\usepackage{color}
\usepackage{algorithm}
\usepackage{algorithmic}
\usepackage{amsthm}
\usepackage{mathrsfs}

\newcommand*\patchAmsMathEnvironmentForLineno[1]{%
  \expandafter\let\csname old#1\expandafter\endcsname\csname #1\endcsname
  \expandafter\let\csname oldend#1\expandafter\endcsname\csname end#1\endcsname
  \renewenvironment{#1}%
     {\linenomath\csname old#1\endcsname}%
     {\csname oldend#1\endcsname\endlinenomath}}%
\newcommand*\patchBothAmsMathEnvironmentsForLineno[1]{%
  \patchAmsMathEnvironmentForLineno{#1}%
  \patchAmsMathEnvironmentForLineno{#1*}}%
\AtBeginDocument{%
\patchBothAmsMathEnvironmentsForLineno{equation}%
\patchBothAmsMathEnvironmentsForLineno{align}%
\patchBothAmsMathEnvironmentsForLineno{flalign}%
\patchBothAmsMathEnvironmentsForLineno{alignat}%
\patchBothAmsMathEnvironmentsForLineno{gather}%
\patchBothAmsMathEnvironmentsForLineno{multline}%
}

\journal{Renewable Energy}

\begin{document}

\begin{frontmatter}
\title{Joint Optimization of Wind Farm Layout Considering Optimal Control}
%
%

\author[label1]{Kaixuan~Chen}
\author[label1]{Jin~Lin\corref{cor1}}
\ead{linjin@tsinghua.edu.cn}
\author[label1]{Yiwei~Qiu}
\author[label1]{Feng Liu}
\author[label1,label3]{Yonghua~Song}

\address[label1]{State Key Laboratory of Control and Simulation of Power Systems and Generation Equipment, Department of Electrical Engineering, Tsinghua University, Beijing 100087, China}
\address[label3]{Department of Electrical and Computer Engineering, University of Macau, Macau 999078, China}

\cortext[cor1]{Corresponding author}
\address{}

\begin{abstract}
  The wake effect is one of the leading causes of energy losses in offshore wind farms (WFs). Both turbine placement and cooperative control can influence the wake interactions inside the WF and thus the overall WF power production. Traditionally, greedy control strategy is assumed in the layout design phase. To exploit the potential synergy between the WF layout and control so that a system-level optimal layout can be obtained with the greatest energy yields, the layout optimization should be performed with cooperative control considerations. For this purpose, a novel two-stage WF layout optimization model is developed in this paper. Cooperative WF control of both turbine yaw and axis-induction are considered in the WF layout design. 
  However, the integration of WF control makes the layout optimization much more complicated and results in a large-scale nonconvex problem, hindering the application of current layout optimization methods. To increase the computational efficiency, we leverage the hierarchy and decomposability of the joint optimization problem and design a decomposition-based hybrid method (DBHM). Instead of manipulating all the variables simultaneously, the joint optimization problem is decomposed into several subproblems and their coordination. Case studies are carried out on different WFs. It is shown that WF layouts with higher energy yields can be obtained by the proposed joint optimization compared to traditional separate layout optimization. Moreover, the computational advantages of the proposed DBHM on the considered joint layout optimization problem are also demonstrated.
\end{abstract}


%
\end{frontmatter}
\section{Introduction}


Wind is one of the fastest growing renewable energy forms. According to the Global Wind Energy Council, the global offshore wind power capacity is expected
to reach 234 GW by the end of 2030  \cite{councilwind}. With this huge amount of installation, the improvement of wind farm (WF) power production has been the central concern for WF development \cite{veers2019grand}. 
The aerodynamic interactions among wind turbines (WT) in a WF, i.e. the wake effect, has been observed to be one of the major causes of offshore WF energy losses \cite{veers2019grand, boersma2017tutorial}. To mitigate the wake effect and increase the energy yields, WF layout design such as placing the turbines farther away in the prevailing wind direction, and cooperative control to coordinate different turbine operation, are two main approaches and have attracted many studies in recent years.

WF control optimization aims to minimize wake losses by cooperatively operating turbines in the WF as a whole. A review of the literature on wake control strategies can be found in \cite{boersma2017tutorial}. Axial induction and wake redirection are the two actuation methods for turbines to control their wakes. In the axial-induction-based control, the generator torques and blade pitch angles are adjusted \cite{ahmadyar2016coordinated}. Alternatively, in the wake redirection control, the upstream turbine is purposely misaligned with the incoming flow with a yaw angle to steer the wake away from the downstream turbines \cite{dar2016windfarm}. Field tests have shown that both methods have potential for increasing the WF power production, though clearly the yaw actuation is much more effective \cite{van2019effects, fleming2017field}. The combination of the two control methods are also investigated \cite{park2016bayesian}. Though WF control methods have been studied quite extensively, most of them are conducted based on regular-shaped WF layouts. The improvement of control optimization for an optimized layout is still questionable.

WF layout optimization considers the stochastic distribution of the overall wind conditions to design a farm layout so that the expected annual energy production (AEP) can be maximized \cite{hou2019review}. Both modeling and optimization algorithms have been studied. Computational Fluid Dynamics (CFD) wake model is accurate, but is too complex to be applied in the WF design process. Analytical parametric wake models are commonly used in the optimization problem due to their low computational costs \cite{GAO2016192} \cite{tao2020wind}.
Various analytical wake models have been employed to characterize wake interactions \cite{tao2019optimal, gonzalez2013new, bastankhah2016experimental}. Refs. \cite{boersma2017tutorial} and \cite{annoni2018analysis} validate the accuracy of different wakes models and provide guidance for the model selection. 
From the optimization modeling perspective, turbine positions can be either modelled by discrete grid points or continuous coordinates \cite{wu2021design}. In the grid model, the area are partitioned into grids and each grid point represents a possible WT position \cite{long2017formulation, cao2020optimizing, turner2014new}. To remove the limitation and get a more accurate layout design model, more works model the turbine locations by continuous coordinates in recent years \cite{MITTAL2016133} \cite{reddy2020wind}.
On the other hand, because the wake model is typically nonlinear and nonconvex \cite{brogna2020new}, difficulties arise in solving these optimization problems. While a few papers have attempted gradient-based methods \cite{PEREZ2013389} \cite{park2015layout}, to overcome the nonconvexity, most studies have employed heuristic optimization methods. Many algorithms have been applied and obtain WF layouts with good AEP improvements such as random search \cite{feng2015solving}\cite{2012Wind}, genetic algorithm \cite{ju2019wind} \cite{huang20183} \cite{YAMANIDOUZISORKHABI2018341}, evolutionary algorithm \cite{wilson2018evolutionary}, ant colony algorithm \cite{EROGLU201253}, particle swarm optimization (PSO)\cite{long2015two, hou2016offshore} and so on. A performance comparison of these methods for layout optimization is provided in \cite{brogna2020new}. 


Generally, most of the above mentioned WF layout design studies \cite{GAO2016192}-\cite{long2015two} do well in mitigating the wake effect and optimizing the WF layouts, but in these studies turbines are assumed to operate individually referred to as greedy operation. With the development of WF cooperative control, the greedy operation assumption may deviate from the WF operation reality in the future \cite{zong2021experimental}. The separation between the cooperative control and layout studies eliminates the possibility of leveraging their potential synergy in the layout design phases. To find system-optimal designs, the layout optimization should considers the cooperative WF control instead of assuming greedy control. 
\cite{hou2016offshore} \cite{wang2016novel} and \cite{pedersen2020integrated} report an attempt in this direction and show the necessity of the joint turbine placement and cooperative control optimization. However, to simplify the optimization, only axis-induction control with restricted farm layouts is considered. The synergy between WF yaw operation and general layout optimization is still questionable. Besides, the computational burden increases dramatically compared to separate layout and control optimizations \cite{fleming2016wind}. 

To address the problem, this paper focuses on the WF layout optimization problem with optimal cooperative control considerations. We advance the current researches by constructing a novel two-stage joint WF layout and control optimization model. Cooperative control for different wind directions in the whole wind rose are integrated into the layout optimization so that the system-level optimal layout can be obtained. However, the joint considerations of WF control and layout make the optimization model a large-scale nonconvex problem. Because the search space grows exponentially with the variable dimension, the computational burden tends to increase infeasible if directly applying existing layout optimization methods to optimize all variables simultaneously. To address this challenge, a decomposition-based hybrid method (DBHM) is proposed to expedite the computation of the joint optimization. Our contribution is twofold: 


1) The WF layout is optimized with optimal cooperative control considerations. A novel two-stage WF optimization model is formulated to jointly optimize the WF control and layout. Both yaw and axis-induction control are considered. The synergy of WF layout and operation can be thoroughly analyzed using the proposed model. Compared to the traditional separate optimizations, increased AEP is achieved.

2) Leveraging the hierarchy and decomposability of the joint optimization problem, a decomposition-based hybrid method is proposed. Instead of manipulating all of the variables simultaneously, the original problem is decomposed into several subproblems and their coordination. Compared to directly applying traditional layout optimization methods in the novel joint optimization problem considered in this paper, the computational burden is alleviated and better results are obtained. 

The remainder of this paper is organized as follows. Section \ref{sec:back} reviews the wake effect and formulates the WF power model. The proposed joint optimization model is constructed in Section \ref{sec:model}. In Section \ref{sec:method}, the proposed DBHM is introduced and applied so that the joint optimization can be solved efficiently. Case studies and conclusions are presented in Sections \ref{sec:case} and \ref{sec:cons}, respectively.

\section{Wake and Wind Farm Models} \label{sec:back}
In this section, the wake model and WF power model are briefly reviewed. 
The Flow Redirection and Induction in Steady-state (FLORIS) model is employed. FLORIS characterises the steady-state wake interaction in the WF and has been widely accepted in WF control and layout studies \cite{hou2019review, doekemeijer2019tutorial, 2016Maximization}. Its accuracy has been extensively verified by both wind tunnel testings \cite{bastankhah2016experimental} and realistic field measurements \cite{annoni2018analysis}. Compared to the classic Jensen model, FLORIS further characterizes the Gaussian shape of the wake and models the yaw actuation.
\begin{figure}[b]
	\centering
		\includegraphics[width=3.0 in]{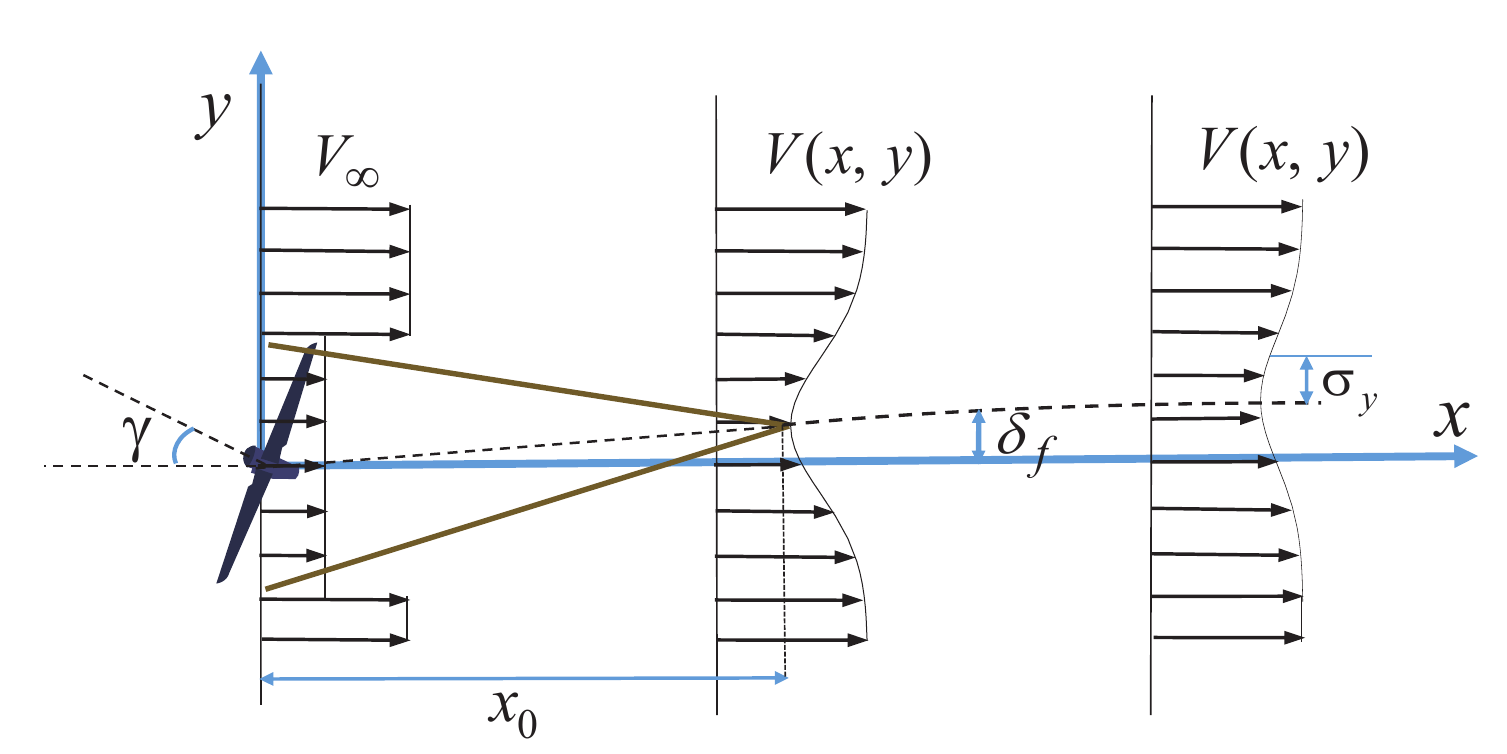}
    \vspace{-7pt}
	\caption{Schematic of the single wake model under yaw condition \cite{bastankhah2016experimental}.}
	\label{fig:wake}
\end{figure}

A schematic overview of a single wake model is given in Fig. \ref{fig:wake} \cite{bastankhah2016experimental}. In this paper, we considers the two dimensional layout optimization, and the FLORIS is utilized. The near-wake zone is modeled as a linearly converging cone with the base at the turbine rotor. In the far wake zone which is the turbine interaction region, the wake deficit follows the shape of a Gaussian distribution:
\begin{align}
  V(x, y)/ V_{\infty} = 1 - (1 - \sqrt{1 - \frac{\sigma_{y0}}{\sigma_y} C_T}) \exp(\frac{(y - \delta_f)^2}{2\sigma_y^2}) \label{wakestart}
\end{align}
\noindent
where $V_{\infty}$ is the freestream wind speed. $x$ and $y$ denote the distance from the rotor center in the wind parallel direction and the lateral direction, respectively. $C_T$ is the thrust coefficient that is related to the turbine operation by an axial induction factor $\alpha$. $C_T = 4\alpha(1-\alpha)$. $\sigma_{y}$ is the standard deviation of the Gaussian shape that is calculated according to
\begin{align}
  \sigma_{y} = \sigma_{y0} + (x - x_0)k_y  ,  \hspace{8pt}  \sigma_{y0} = \frac{D}{2\sqrt{2}} \cos \gamma,
\end{align}
\noindent
where $k_y$ is the wake expansion coefficient. $x_0$ is the length of the near-wake zone \cite{doekemeijer2019tutorial}. $\gamma$ is the yaw angle of the turbine. $D$ is the rotor diameter. Due to the yaw misalignment, the Gaussian-shaped wake centerline $\delta_f$ is deflected from the rotor center, and is modeled as
\begin{align}
  \delta_f = a_d \cdot D + b_d \cdot x  + \tan(\phi) x_0 + \frac{\phi}{5.2}E_0 \times
  \sqrt{\frac{\sigma_{y0}}{k_y C_T}} \cdot \ln[ \frac{(1.6 + \sqrt{C_T})(1.6\sqrt{\sigma_{y}/\sigma_{y0}} - \sqrt{C_T})}{(1.6 - \sqrt{C_T})(1.6\sqrt{\sigma_{y}/\sigma_{y0}} + \sqrt{C_T})}] \label{wakem}
\end{align}
\noindent
where $E_0 = C_0^2 - 3e^{1/12}C_0+3e^{1/3}$ and $C_0 = 1-\sqrt{1-C_T}$. $e$ is the natural constant. $a_d$ and $b_d$ are tuning parameters. $\phi$ is a function of $\gamma$ that is given by:
\begin{align}
  \phi \approx \frac{0.3 \gamma}{\cos \gamma} (1- \sqrt{1 - C_T \cos \gamma})
\end{align}
\noindent
\noindent
Equation (\ref{wakem}) gives the value of the wake deflection at downwind locations as a function of turbine control variable $\gamma$ and $C_T$ as well as the wake empirical parameters. Equation (\ref{wakem}) takes into account both the wake deflection caused by the yaw angle, and the deflection caused by the rotor rotation. Readers interested in the formulation details are referred to Appendix B in \cite{bastankhah2016experimental}.
The combination of multiple wakes is calculated by the traditional sum of squares method \cite{park2015layout}. $k_y$, $a_d$ and $b_d$ are the tuning parameters of this wake model for different WFs. 

The power generated by each turbine under yaw conditions can be computed based on the actuator disk theory \cite{annoni2018analysis}.
\begin{align}
  P_{\textrm{WF}} = \sum_{i=1}^{N} P_{i}= \sum_{i=1}^{N} \frac{1}{2} \rho (\frac{\pi}{4} D^2) C_{P_i} \cos (\gamma_i)^{p_p} \overline{V}_{i}^{3}   \label{wakeend} 
\end{align}
\noindent

$P_{\textrm{WF}}$ is the WF power production calculated by the summation of each turbine power $P_i$. $N$ is the number of turbines in the WF and the subscript $i$ is the index for the $i$th WT. $\rho$ is the air density. $\overline{V}_{i}$ is the velocity  $V(x, y)$ averaged over the rotor disk. $p_p$ is the tunable parameter for turbine power production under yaw conditions. $C_{P}$ is the power coefficient that is also directly related to the axial induction factor $\alpha$: 
\begin{align}
 C_{P_i} = 4\alpha_i(1-\alpha_i)^2 \label{cp}
\end{align}
\noindent

As commonly employed in WF studies, $\gamma$ and $\alpha$ of each turbine are considered as the control variables and are used to regulate the WF operation.
The upstream turbine influences the wind speed faced by downstream turbines, and the WF is thus operated as a whole.

Through the FLORIS Gaussian modeling, the wake interaction and the WF power production are modeled as a smooth function in terms of turbine locations \cite{park2015layout, doekemeijer2019tutorial}.
Readers interested in more details of the FLORIS wake model and its verification are referred to \cite{bastankhah2016experimental, annoni2018analysis}.

\section{WF Layout Optimization Model with Cooperative Control Considerations} \label{sec:model}
Traditionally, WF layout and cooperative control optimization are considered separately and sequentially. In traditional layout optimization models, the ``greedy" WT operation settings are assumed, i.e. $\gamma = 0^{\textrm{o}}$ and $\alpha = 1/3$ for all turbines \cite{hou2019review}\cite{brogna2020new}. With the development of WF cooperative control strategy, the greedy operation assumption may deviate from the WF operation reality in the future. Thus the layout optimization results may not find system-level optimal designs, because the couplings of system design and operation are ignored in the design phase \cite{fathy2001coupling}. To exploit the potential synergy, a two-stage optimization model is constructed so that the WF layout optimization can consider the cooperative control strategy. 
\subsection{Two-Stage WF Layout Optimization with Operation Considerations}
In this paper, we focus on the design purpose of maximizing the potential annual energy production (AEP) capacity of the whole WF. For this purpose, the probability distribution of the stochastic freestream wind direction on the target site should be considered. The first-stage WF layout optimization is formulated as follows:
\begin{align}
  &\min \limits_{\textbf{x}, \textbf{y}} \hspace{8pt}  -f_{\textrm{AEP}}(\textbf{x}, \textbf{y})  
   = - T \sum_{\omega = 1}^{N_\theta} p(\theta_{\omega}) P_{\textrm{WF}}^{\ast} (\textbf{x}, \textbf{y}, \theta_{\omega}) \label{obj} \\
  &\textrm{s.t.:} \hspace{7pt} (x_i - x_j)^2 + (y_i - y_j)^2 \geqslant L^2, \hspace{8pt}  i,j =1, \cdots, N  \label{poscons} \\
   &\hspace{55pt} x_l \leq x_i \leq x_u, \hspace{8pt} i =1, \cdots, N \label{terrcons} \\
   &\hspace{55pt} y_l \leq y_i \leq y_u, \hspace{8pt} i =1, \cdots, N \label{terrcons2}    
\end{align}
\noindent

$f_{\textrm{AEP}}$ is the WF AEP, which is modelled by the product of the expected power production and the total hours in one year \cite{tao2019optimal} \cite{2016Maximization}. $N_{\theta}$ is the total number of intervals with equal width into which the wind direction is discretized, and $\omega$ is the wind direction index. $\theta_{\omega}$ denotes the wind direction for scenario $\omega$. $p(\theta_{\omega})$ is the frequency of occurrence for the $\omega$-th wind direction. The probability distribution of the wind direction is typically recorded by annual measurements \cite{tao2019optimal} \cite{park2015layout}. $T=8760 \textrm{h}$ is the number of hours in one year. Because the freestream wind speed has little impact on the WF wake optimization \cite{2016Maximization} \cite{zhong2016decentralized}, we do not consider the speed distribution. The decision variables $\textbf{x} = [x_1, x_2, ..., x_{N}]$, $\textbf{y} = [y_1, y_2, ..., y_{N}]$ represent the turbine coordinates. $N$ is the total number of turbines in the WF, and subscripts $i$ and $j$ are the WT indexes. 
(\ref{poscons}) requires that the distance between every two WTs should be larger than a specified distance $L$ for the safe operation \cite{tao2019optimal}. The terrain boundaries are constrained in (\ref{terrcons}) and (\ref{terrcons2}), where the subscripts $l$ and $u$ denote the lower and upper boundaries respectively. While a rectangular terrain is modeled for convenience in this paper, other geometric shapes can be easily considered \cite{reddy2021efficient}. 

$P^{\ast}_{\textrm{WF}}$ denotes the optimal WF power production for a given $\textbf{x}, \textbf{y}$ and wind scenario $\theta_{\omega}$, which is the optimization result of the second-stage WF control optimization. For the $\omega$th inflow wind scenario, the control optimization is given by:
\begin{align}
  P^{\ast}_{\textrm{WF}}(\textbf{x}, \textbf{y},  \theta_{\omega})=  \max \limits_{\bm{\gamma}_{\omega}, \bm{\alpha}_{\omega}}  \sum_{i=1}^{N} \frac{1}{2} \rho (\frac{\pi}{4} D^2) C_{P_i} \cos (\gamma^{\omega}_i)^{p_p} \overline{V}_{i}^{3} \label{controlobj}\\
   \textrm{s.t.:} \hspace{8pt} \gamma_{min} \leq  \gamma^{\omega}_i  \leq \gamma_{max}, \hspace{20pt}  \label{opecons} i =1, \cdots, N \hspace{40pt} \\
     \alpha_{min} \leq  \alpha^{\omega}_i  \leq \alpha_{max}, \hspace{20pt}  \label{opecons2} i =1, \cdots, N \hspace{40pt}\\
    \textrm{wake interaction}: (\ref{wakestart})-(\ref{cp})  \hspace{60pt} \label{wakeconstraint}
\end{align}

The yaw angles $\bm{\gamma}_{\omega} = [\gamma_{1}^{\omega}, ..., \gamma_{i}^{\omega}, ..., \gamma_{N}^{\omega}]$ and axial induction factors $\bm{\alpha}_{\omega} = [\alpha_{1}^{\omega}, ..., \alpha_{i}^{\omega}, ..., \alpha_{N}^{\omega}]$ are the control variables of all turbines for a given $\textbf{x}, \textbf{y}$ under wind scenario $\theta_{\omega}$. Considering the wake interactions inside the WF in (\ref{wakeconstraint}), (\ref{controlobj}) tries to maximize the power production by coordinate different turbine operation. Constraints (\ref{opecons}) and (\ref{opecons2}) limit the lower and upper bounds of the control actions as commonly modelled by WF control researches \cite{boersma2017tutorial}. 

\begin{figure}[htb]
	\centering
  \vspace{-7pt}
		\includegraphics[width=2.5 in]{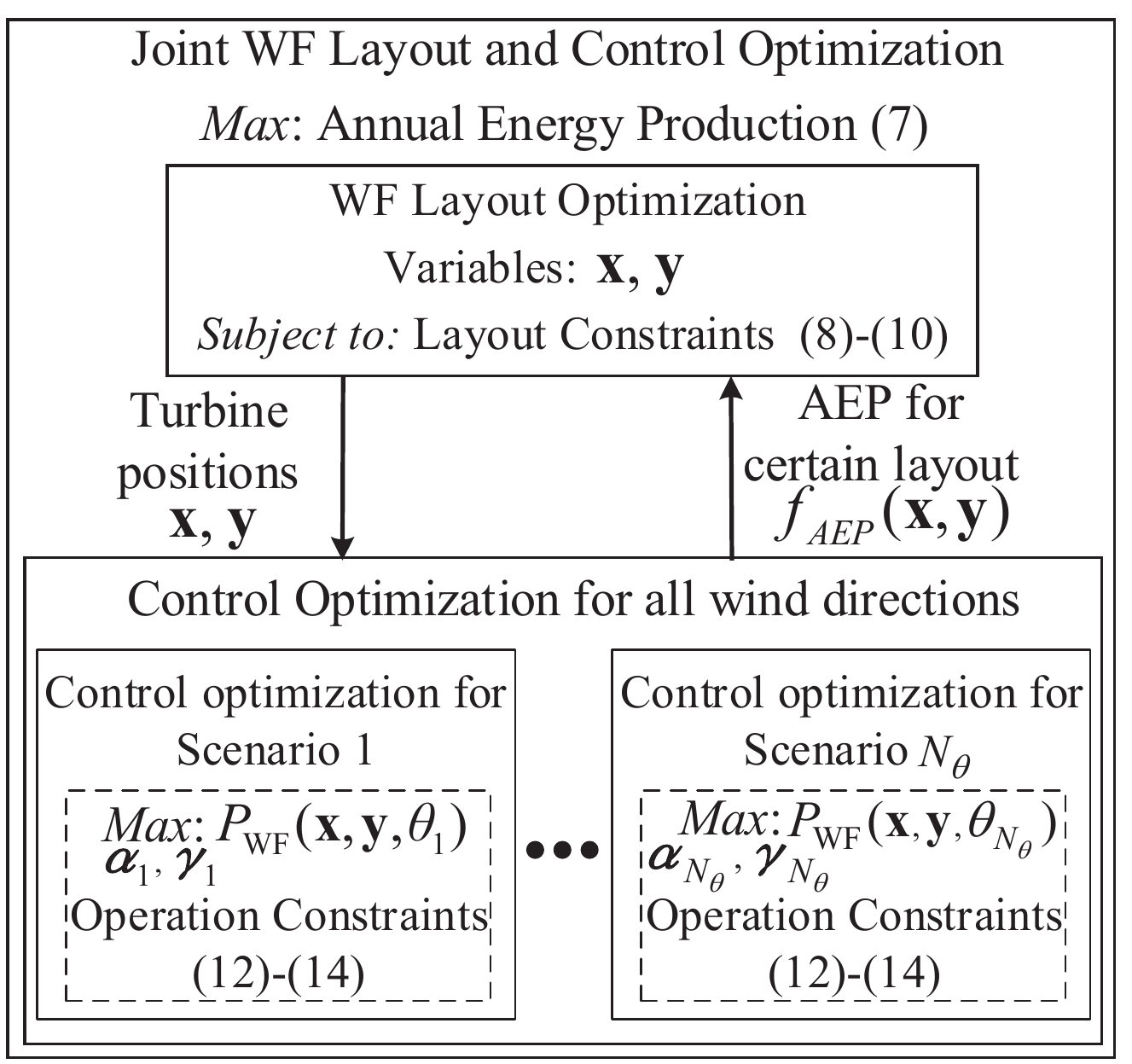}
  \vspace{-5pt}
	\caption{Framework of the joint layout and control optimization model.}
	\label{fig:modelstructure}\vspace{-12pt}
\end{figure}

(\ref{obj})-(\ref{wakeconstraint}) together formulate the joint optimization model that can simultaneously considers the WF layout and operation. The hierarchy diagram of the two-stage optimization problem is shown in Fig. \ref{fig:modelstructure}. The first-stage layout problem (\ref{obj})-(\ref{terrcons2}) manipulates turbine coordinates. The second-stage operation problems (\ref{controlobj})-(\ref{wakeconstraint}), which can be viewed as constraints of the first-stage layout optimization, determine the optimal control for each wind direction considering wake interactions among WTs. To account for all wind directions in the layout design in (\ref{obj}), for a given $[\textbf{x}, \textbf{y}]$, $N_{\theta}$ control optimizations are required to be solved to evaluate $f_{\textrm{AEP}}$. Thus, the total number of the decision variables of the joint model is $(N_{\theta} + 1)*2N$ with $2N $ layout variables and $N_{\theta} *2N$ operation variables.

\subsection{Complexity analysis of the Joint Optimization Problem} \label{comptime}
Because of the wake interactions (\ref{wakestart})-(\ref{wakeend}), the WF optimization model is nonlinear and nonconvex \cite{hou2019review}\cite{brogna2020new}.
For traditional isolated layout or control optimizations that considers only $2N$ variables, heuristic algorithms are employed by most studies due to
their global exploration ability \cite{hou2019review} \cite{brogna2020new} \cite{ju2019wind}. However, heuristic methods are known to suffer from slow convergence, and the computational burden tends to increase exponentially with the expansion of the problem size \cite{park2015layout}\cite{2010Runtime}, making the WF layout optimization typically time-consuming \cite{brogna2020new} \cite{pedersen2020integrated}. 

For the two-stage joint optimization model (\ref{obj})-(\ref{wakeconstraint}), in each layout evaluation step, the second-stage operation problem is required to optimize WF control over the entire distribution of wind directions. Such integration of the scenario-specific $N_{\theta} *2N$ control variables multiplies the problem size of the already time-consuming layout design, making the joint optimization model a much larger nonconvex problem. Thus, directly applying traditional layout optimization methods manipulating all variables simultaneously is computational expensive for the two-stage optimization \cite{2010Runtime}\cite{pedersen2020integrated}, as shown in the experimental results of Section \ref{sec:case}.

On the other hand, as shown in Fig. \ref{fig:modelstructure}, though the joint optimization model is complex considering the entire wind direction distribution, the operation variables $\bm{\gamma}_{\omega}, \bm{\alpha}_{\omega}$ and constraints (\ref{opecons})-(\ref{opecons2}) for different wind scenarios are not coupled. When the first-stage layout decisions, $\textbf{x}$ and $\textbf{y}$, are given, the second-stage operation problem naturally decouples into low-scale subproblems for each wind scenario which can be efficiently solved in parallel.
This suggests that a structuring of the solving algorithm to take into account the hierarchy and decomposability of the problem is a promising direction.


\vspace{-3pt}

\section{Decomposition Formulation of the Joint Optimization Model} \label{sec:method}
In this section, we propose a decomposition-based hybrid method (DBHM) so that the joint WF layout and operation optimization can be solved efficiently.
The proposed DBHM is composed of two steps that combine the heuristic exploration and the problem decomposition.

First, the heuristic optimization method is employed to optimize the isolated WF layout, the results of which are used to ``warm start" the joint optimization in the next step. Heuristic approaches are adopted because they do not get trapped at local optima arising from the nonconvexity. In this step, PSO is adopted which has been wildly used by isolated WF layout studies\cite{hou2015optimized, hou2016offshore}. A penalty function is adopted to handle the layout distance constraints \cite{POOKPUNT2013266}. The penalty factor, PF, is $10^{5}$ and is adjusted by trial and error. 
\begin{align}
  \min  \limits_{\textbf{x}, \textbf{y}} \hspace{8pt} -f_{\textrm{AEP}} -\textrm{PF} \sum_{i=1}^{N}\sum_{j=i}^{N} \min((x_i - x_j)^2  + (y_i - y_j)^2 - (L)^2, 0)
\end{align}

Furthermore, a decomposition formulation of the joint optimization model (\ref{obj})-(\ref{wakeconstraint}) is constructed. The original problem (\ref{obj})-(\ref{wakeconstraint}) is decomposed into a set of subproblems that can be solved efficiently through an iterative process of subproblem optimizations and the coordination among them. 

For the original joint optimization model, the only couplings among different wind scenarios are the turbine coordinates $\textbf{x}$ and $\textbf{y}$. These coupled variables are handled following the idea of dual decomposition \cite{boyd2011distributed}. First, local copies of the turbine coordinates, $\textbf{x}^{\omega} = [x_1^{\omega}, x_2^{\omega}, ..., x_{N}^{\omega}]$ and $\textbf{y}^{\omega} = [y_1^{\omega}, y_2^{\omega}, ..., y_{N}^{\omega}]$ are introduced to each wind scenario $\omega$ together with a set of consistency constraints. In this way, the original joint optimization (\ref{obj})-(\ref{wakeconstraint}) can be written equivalently as follows: 
\begin{align}
   & \min  \hspace{4pt}  - T \sum_{\omega = 1}^{N_\theta}  P_{\textrm{WF}} (\textbf{x}^{\omega},\textbf{y}^{\omega}, \bm{\gamma}_{\omega}, \bm{\alpha}_{\omega}, \theta_{\omega})  p(\theta_{\omega}) \label{admobj} \\
   & \hspace{4pt}  \textrm{s.t.:} \hspace{12pt} x_i^{\omega} = x_i, \hspace{7pt} y_i^{\omega} = y_i, \hspace{7pt} \forall i,\omega \label{admmconsesus} \\
   & \hspace{30pt} \textrm{Layout Constraits}: (\ref{poscons}) - (\ref{terrcons2}) \label{admdis} \\
   & \hspace{30pt} \textrm{Operation Constraits}: (\ref{opecons}) - (\ref{wakeconstraint}), \hspace{4pt} \forall \omega \label{admmope}
\end{align}

The consensus constraint (\ref{admmconsesus}) enforces that the introduced new variables $\textbf{x}^{\omega}$ and $\textbf{y}^{\omega}$ of each scenario are equal to the original layout variables $\textbf{x}$ and $\textbf{y}$. The couplings among different wind scenarios are transformed from layout variables to the consensus constraint (\ref{admmconsesus}).


Then, the method of multipliers is adopted to relax the consensus constraints in the form of augmented Lagrangian function. In this way, (\ref{admobj})-(\ref{admmope}) can be decomposed into a two-block formulation. One block is a coordination problem updating $\textbf{x}$ and $\textbf{y}$. The second block considers $N_{\theta}$ individual scenario subproblems.
For the $\omega$th inflow wind scenario, the corresponding subproblem can be constructed as follows  
\begin{align}
  & \min \limits_{\varepsilon^{\omega},\textbf{x}^{\omega},\textbf{y}^{\omega}, \bm{\gamma}_{\omega}, \bm{\alpha}_{\omega}}  -P_{\textrm{WF}} (\textbf{x}^{\omega},\textbf{y}^{\omega}, \bm{\gamma}_{\omega}, \bm{\alpha}_{\omega}, \theta_{\omega})  p(\theta_{\omega}) + \varepsilon^{\omega}  \label{subobj}\\
  & \hspace{15pt} \textrm{s.t.:} \hspace{25pt} (x_i^{\omega} - x_j^{\omega})^2 + (y_i^{\omega} - y_j^{\omega})^2 \geqslant (L)^{2}, \hspace{5pt} \forall i,j   \label{addis}\\
  & \hspace{90pt} x_l \leq x_i^{\omega} \leq x_u, \hspace{8pt}\forall i \label{adup}\\
  & \hspace{90pt} y_l \leq y_i^{\omega} \leq y_u, \hspace{8pt}\forall i  \label{adlow} \\
  & \hspace{70pt} \textrm{Operation Constraits}: (\ref{opecons}) - (\ref{wakeconstraint}) \label{adop}\\
  & \hspace{5pt} \varepsilon^{\omega} \geq \sum_{i=1}^{N} \lambda_{i,x}^{\omega}(x_i - x_i^{\omega}) + \lambda_{i,y}^{\omega}(y_i - y_i^{\omega})  + \mu (x_i - x_i^{\omega})^2+ \mu(y_i - y_i^{\omega})^2  \label{multip}
\end{align}

The decision variables of the subproblem are all local variables $(\textbf{x}^{\omega},\textbf{y}^{\omega}, \bm{\gamma}_{\omega}, \bm{\alpha}_{\omega})$ that are independent of other wind scenarios. Constraints (\ref{addis})-(\ref{adlow}) limit the feasible region of the introduced local variables $\textbf{x}^{\omega}$ and $\textbf{y}^{\omega}$ in the same manner as  (\ref{poscons})-(\ref{terrcons2}). Because in the proposed solving framework the variables are updated through an iterative process (described in Subsection \ref{DBHM}), original linking variables $\textbf{x}$ and $\textbf{y}$ in (\ref{multip}) are parameters in the subproblem. $\varepsilon^{\omega}$ is an ancillary variable that formulate the augmented Lagrangian function with (\ref{multip}) to relax the consensus constraint (\ref{admmconsesus}). The violations of (\ref{admmconsesus}) are penalized in the objective function by the augmented Lagrangian form \cite{boyd2011distributed}: $\lambda_{i,x}^{\omega}$ and $\lambda_{i,y}^{\omega}$ are the Lagrange multipliers corresponding to the consensus constraints in (\ref{admmconsesus}). $\mu$ are the corresponding penalty factors. For each subproblem, only one particular wind condition is considered and the number of decision variables is only $4N$. Because the subproblems (\ref{subobj})-(\ref{multip}) for different wind scenarios are mutually independent, the second block naturally decomposes into $N_{\theta}$ separate subproblems. 

On the other hand, the formulation of the first-block coordination problem is given as 
\begin{align}
  \min \limits_{\textbf{x}, \textbf{y}}  \sum_{\omega = 1}^{N_\theta}  \sum_{i = 1}^{N} \lambda_{i,x}^{\omega}(x_i - x_i^{\omega}) + \lambda_{i,y}^{\omega}(y_i - y_i^{\omega}) + \mu (x_i - x_i^{\omega})^2+ \mu (y_i - y_i^{\omega})^2  \label{master}
\end{align}
\noindent

The coordination problem updates the original layout variable $\textbf{x}$ and $\textbf{y}$ based on the subproblem optimization results (\ref{subobj})-(\ref{multip}) of all wind scenarios. The coordination problem (\ref{master}) is a convex quadratic programming with no constraints that can be solved analytically.
\begin{align}
   x_i = \frac{1}{N_{\theta}} \sum_{\omega = 1}^{N_{\theta}} (x_i^{\omega} - \lambda_{i,x}^{\omega} /2\mu  ), \hspace{8pt}  \forall i \label{upperx}\\
   y_i = \frac{1}{N_{\theta}} \sum_{\omega = 1}^{N_{\theta}} (y_i^{\omega} - \lambda_{i,y}^{\omega} /2\mu ), \hspace{8pt} \forall i \label{uppery}
\end{align}

The large-scale joint optimization problem (\ref{obj})-(\ref{wakeconstraint}) are decomposed into $N_{\theta}$ low-scale individual scenario subproblems (\ref{subobj})-(\ref{multip}) and a simple coordination problem (\ref{master}). The subproblems are mutually independent and can thus be efficiently solved in parallel or sequentially. Through the decomposition formulation,  the two-stage joint optimization problem (\ref{obj})-(\ref{wakeconstraint}) is successfully decomposed into multiple subproblems, each of which considers only one particular wind condition. The WF control optimization can be efficiently considered in the layout design phase to exploit their potential synergy. The detailed iterative solving process of the decomposition problem are described in Appendix A along with the convergence analysis for interested readers.

\section{Case Study} \label{sec:case}


In this section, the synergy of WF operation and layout optimization is analyzed by comparing the joint optimization results with separate layout and control optimizations. 
The simulation is based on the open-source FLORISSE-M platform\cite{doekemeijer2019tutorial}, which has been widely used in WF control and layout studies  \cite{hou2019review} \cite{doekemeijer2019tutorial} \cite{2016Maximization}. Its accuracy has been validated previously by both wind tunnel testings \cite{bastankhah2016experimental} and realistic field measurements \cite{annoni2018analysis}.
The commonly-used NREL Type III WT are assumed to be installed\cite{jonkman2009definition}. The main wake and turbine parameters are provided in Table \ref{Tab1}. For the research
purpose, the yaw angles are allowed to change between $-30^o$ to $30^o$ according to previous researches and realistic experiments \cite{park2016bayesian} \cite{campagnolo2020wind}. Larger yaw range may cause heavier loads to turbines. 

\begin{table}[h]
\centering
\caption{Main Wake and Turbine Parameters}
\label{Tab1}
\begin{tabular}{c c c c c}
\hline
  Air density $\rho$ & $D$ & $p_p$ &  $\alpha_{min}$  & $\alpha_{max}$  \\
  1.29 kg/$m^{3}$ & 126 m & 1.88 &   0.1 & 1/3  \\
\hline
  $a_d$ & $b_d$ & $k_y$  & $\gamma_{min}$ & $\gamma_{max}$ \\
  -0.0356 &  -0.01 & 0.0229 &   $-30^{o}$ & $30^{o}$ \\
\hline
\end{tabular}
\end{table}

\subsection{Illustrative Case}
\begin{figure}[h]
	\centering
  \vspace{-5pt}
		\includegraphics[width=4.2 in]{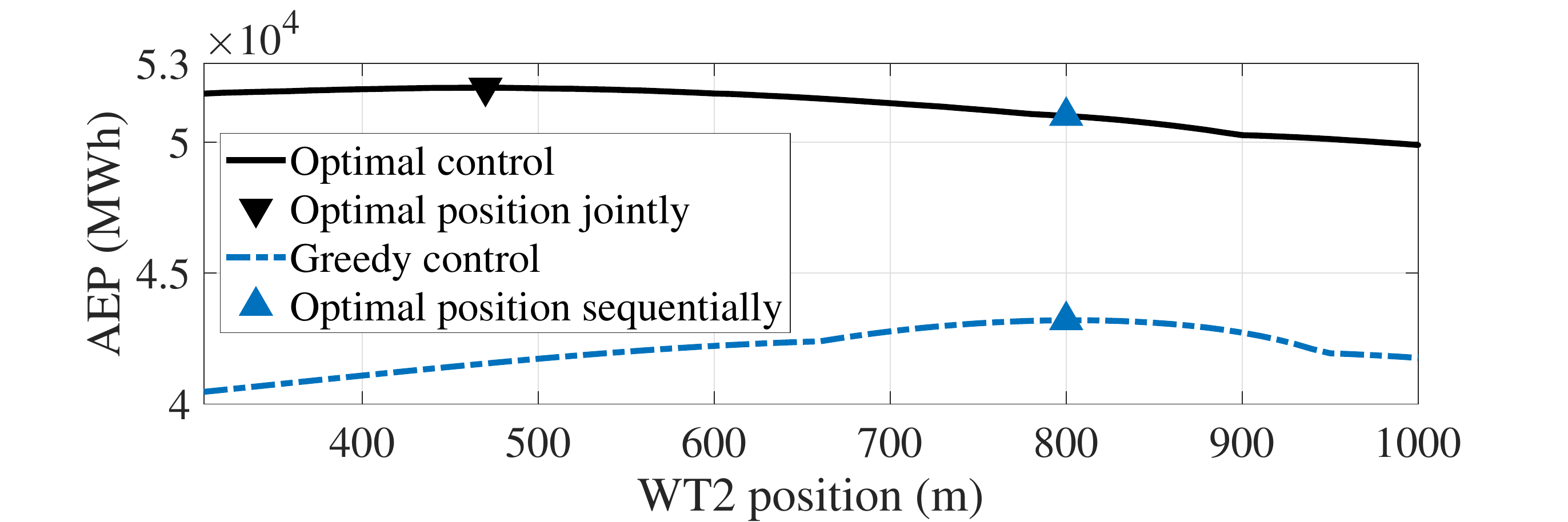}
  \vspace{-4pt}
	\caption{WF potential energy production with respect to WT2 positions with/without cooperative control considerations.}
 \label{fig:3wtillustration}
\end{figure}

To visualize the necessity of considering the optimal operation in the layout design phase, we take an example consisting of 3 WTs on a row. WT1 and WT3 are spaced by 1100~m. The AEP of the three WTs is plotted in Fig. \ref{fig:3wtillustration} with respect to the varying WT2 positions with a 9 m$/$s inflow aligned with the row. This simplified case can be analyzed via 'brute force'. The results at each location are exhausted and visualized.

In Fig. \ref{fig:3wtillustration}, WT2 positions with/without optimal control considerations are compared. The dashed line indicates the results when all of the turbines are operated greedily, which is the common assumption in traditional isolated layout optimizations \cite{hou2019review} \cite{park2015layout}. The optimal position of WT2 in this case is at 800~m, as highlighted by the blue triangle. The corresponding AEP is $4.32 \times 10 ^{4}$~MWh. If operations are then optimized on this turbine position, the AEP is $5.10 \times 10 ^{4}$~MWh.
The black line shows the AEP for optimized WF operation. The optimal position of WT2 in this case is at 470~m as highlighted by the black inverted triangle. The AEP at this position is $5.21 \times 10 ^{4}$~MWh. The optimal yaw angle of the three turbines from WT1 to WT3 in this case are $30^o$, $5.2^o$ and $0^o$, respectively. The optimal axis induction factors are $0.33$, $0.30$ and $0.33$ respectively, which is the same as in greedy control. Compared to the separate layout and operation considerations, the AEP is improved by $2.06 \%$. 

If the minimal turbine distance $L = 4D$ is considered \cite{brogna2020new}\cite{hou2016offshore}, the eligible range for the WT2 position is between 504m and 596m. If the operation and layout are considered sequentially and separately, the optimal position is at the right endpoint (i.e. 596 m) and the corresponding AEP is $5.185 \times 10 ^{4}$~MWh. If operations are considered in the layout design phase, the optimal position of WT2 is at the opposite endpoint (i.e. $504$~m) and the AEP is $5.205 \times 10^{4}$~MWh. The AEP is improved by $0.39\%$.

Comparing the two cases, the WF operation has a clear impact on the optimal turbine position result. To consider such synergy, turbine placement and cooperative control should be simultaneously optimized so that system-level optimal layout can be obtained. 

\subsection{Optimization Results and Analysis} \label{16wtcase}
\begin{figure}[t]
  \vspace{-5pt}
	\centering
		\includegraphics[width=2.6 in]{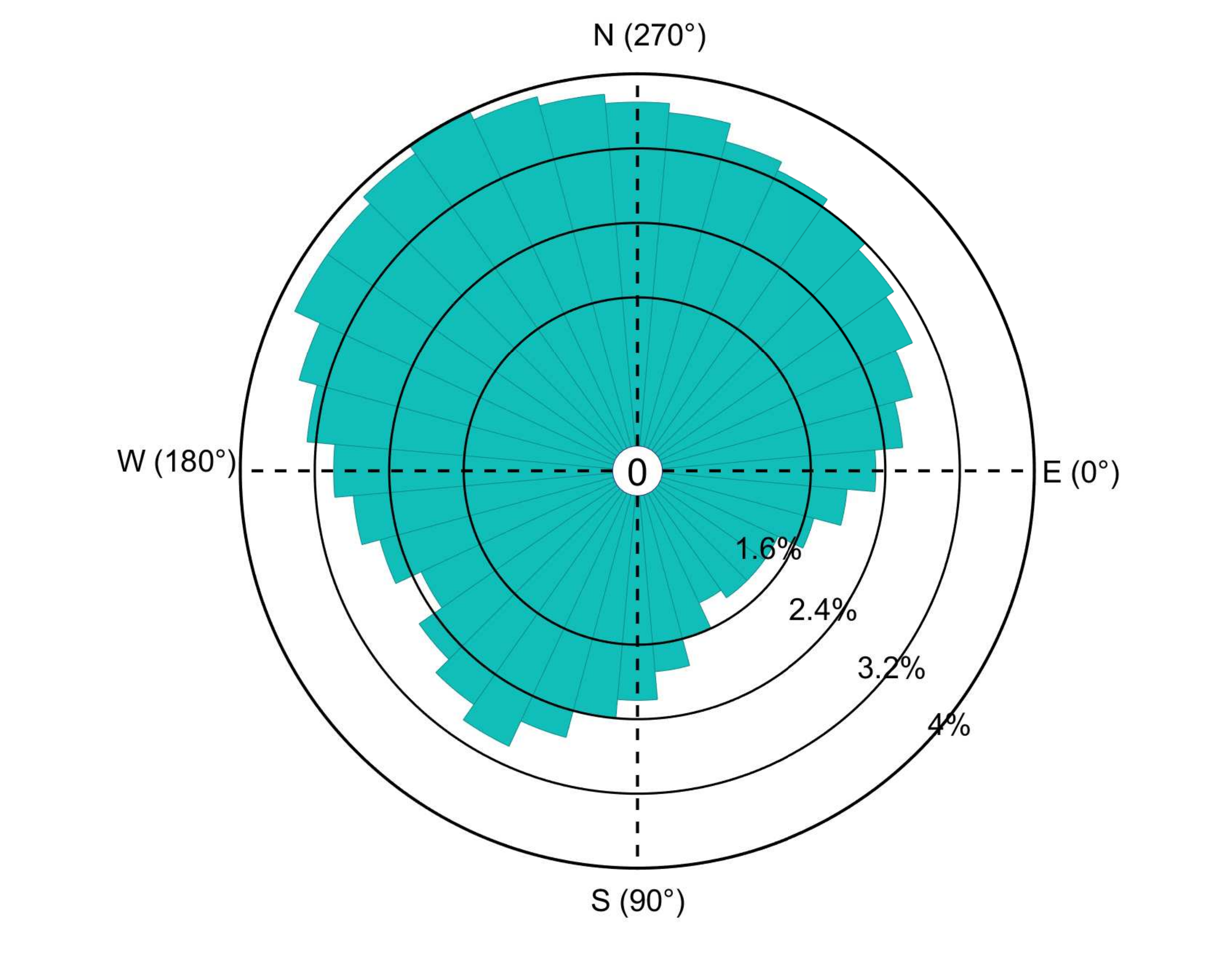}
    \vspace{-13pt}
	\caption{Wind direction probability distribution.}
  \vspace{-10pt}	
\label{fig:windrose}
\end{figure}

A 16-turbine WF located in the offshore area to the west of Denmark is first studied. The wind direction distribution is illustrated in Fig. \ref{fig:windrose} in the wind rose format generated based on the realistic measurements \cite{park2015layout}. The whole range is discretized into 36 wind direction scenarios. The wind speed is assumed to be 9 m/s. The minimal turbine distance $L = 4D$ \cite{brogna2020new}\cite{hou2016offshore}. The land for planning has a square shape with an area of $1900$m $\times$ $1700$m. The initial layout of WF is settled in a rectangular layout that is typically employed in WF studies \cite{ahmadyar2016coordinated} \cite{dar2016windfarm}.

To investigate the synergy of the WF layout and operation, 5 optimization cases are conducted and compared. As the base case, the initial layout with greedy control is studied. Then, in cases 2 and 3, the WF layout and operation are optimized, respectively. Case 4 conducts optimizations sequentially without considering the potential synergy. First, the topology is optimized assuming the greedy WT operation the same as case 3. Then, control is optimized based on the optimized layout.
For the isolated layout optimization, PSO that has been adopted by various WF layout studies is utilized \cite{hou2015optimized, hou2016offshore}. Ten trials are conducted by PSO and the best result is chosen.
The proposed DBHM is finally applied in case 5, where the joint layout and operation optimization (\ref{obj})-(\ref{wakeconstraint}) is performed. 
The results of all of the investigations are summarized in Table \ref{Tab2}.

\begin{table}[h]
\centering
\caption{Results obtained by Various Optimization Considerations}
\label{Tab2}
\begin{tabular}{c c c c c}
\hline
  Case & Layout & Operation & AEP(GWh) & Improvement \\
\hline
  1 & Initial & Greedy & 366.52 &   -- \\
  2 & Initial & Optimized &   373.96 &  2.03\% \\
  3 & Optimized & Greedy &   376.24 &   2.65\% \\
  4 & Optimized & Optimized (Sequent) &   386.49 &   5.45\% \\
  5 & Optimized & Optimized (Joint) &   402.96 &  9.94\% \\    
\hline
\end{tabular}
\end{table}
Comparing case 2 and 3 with case 1, it is found that the isolated layout or control optimization can both increase the WF power production capacity. For optimized WF layout of case 3, the  cooperative control optimization is still profitable. The AEP of case 4 is increased by $ 2.72\%$ compared to case 3.
The gain of considering the optimal control in the layout design phase can be verified by comparing the results of cases 4 and 5.
By exploiting the synergy, better WF layouts  with higher power production are obtained in case 5. Compared to the sequential and separate optimization, the AEP is increased by $4.26\%$ by considering layout and control jointly. For the optimized control results, the upstream turbines purposely yaw to deflect the wake so that it will not at all or partially overlap the downwind turbines and more power can be harvested. On the other hand, the optimal axis induction factors are all close to 1/3, which is the same as in the greedy control. Therefore, for the cooperative WF control, the AEP increase mainly comes from the non-zero yaw angles. This results agree well with the previous WF control studies \cite{E2017Statistical} \cite{munters2018dynamic} and experiments \cite{2016Wind}. Therefore, for the joint control and layout optimization in the paper, the axial induction factors can be fixed at 1/3. Only the synergy between the yaw control and the layout is required to be considered for AEP maximization.

\begin{figure}[h]
	\centering
		\includegraphics[width=5.5 in]{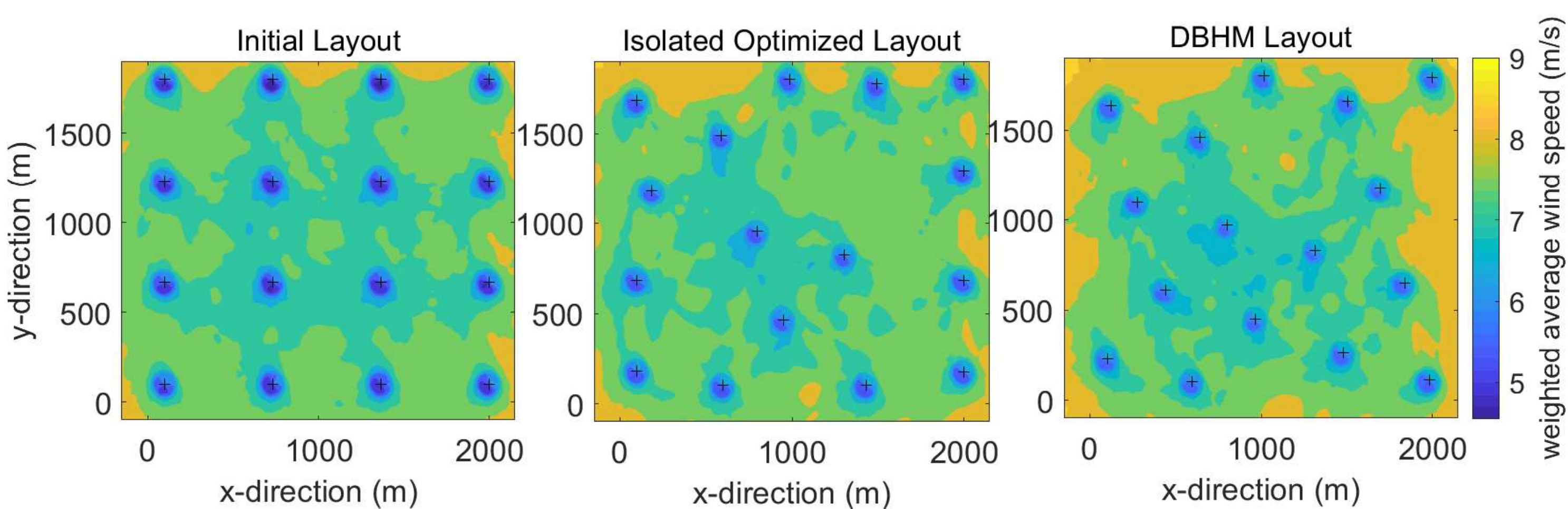}	
  \caption{WF layouts and weighted average wind velocity distributions inside the WF for different cases. It can be shown that the averaged speed on the turbine positions is increased by the DBHM layout.}
	\label{fig:comparelayout}	
\end{figure}

Fig. \ref{fig:comparelayout} depicts the WF layout and the weighted average speed inside the WF for different cases. The speed inside the WF is computed for each of the 36 directions based on the layout and control strategy of corresponding cases. Then, the speed distribution inside the WF is averaged weighted on the probability of the 36 directions in Fig. \ref{fig:windrose}. The initial layout is the typical regular layout. The optimized layout is irregular and scattered. The scattered turbines prevent the formation of a long wake chain and the downstream turbines thus experience less intensified wake. Besides, compared to the isolated optimized layout, the joint optimization result is more scattered. It looks like that the turbines tend to spread in the NNW-SSE and SW-NE directions. This allows for more space between the turbines and better wake recovery. It can be seen that the averaged speed on the turbine positions is increased by the proposed joint optimization and an increased AEP is thus obtained.

The computational efficiency of DBHM is further demonstrated. For the isolated WF layout optimization by PSO in case 3, the average computational time of ten trials is  14443.98~s. For the joint optimization (\ref{obj})-(\ref{wakeconstraint}), Table \ref{Tab3} compares the computation performance of different methods. All of the cases are solved in a MATLAB R2020a environment with a 4 core Intel i7 CPU of 3.6 GHz. Commonly-employed methods for isolated layout optimization, the heuristic PSO \cite{brogna2020new} and the gradient-based sequential convex programming (SCP) \cite{park2015layout}, are applied to compare the efficiency on the considered joint optimization problem. The two methods are directly applied to the joint optimization model (\ref{obj})-(\ref{wakeconstraint}). Thus, the total number of decision variables is 1184. 
The same initial points from the results of case 4 and the terminal condition $\epsilon_{\textrm{tol}} = 10$ are set for SCP and DBMM. For PSO, ten runs are conducted with Matlab Global Optimization toolbox \cite{2010Global}, and the average result is presented because of its uncertainty. 
Though performing well for the WF layout optimizations, the computational performance deteriorates for the joint optimization problem.

\begin{table}[h]
\centering
\caption{Computational Performance Comparison of Different Optimization Methods}
\label{Tab3}
\begin{tabular}{c c c c}
\hline
  Method & SCP & PSO & DBHM \\
\hline
  Average Computation Time(s) & 16006.16 & 94064.74 &  12969.60  \\
  AEP(GWh) & 391.31 & 377.21 &   402.96  \\
  Improvement compared to case 4 & 1.24\% & -2.40\% &  4.26 \% \\
\hline
\end{tabular}
\end{table}

For the DBHM, 36 subproblems are decomposed. $\mu = 10$. For each subproblem, only 64 variables need to be manipulated. 
In this paper, the subproblems are solved in sequence. A parallel implementation will generally lead to more computational savings. 
Because of the nonconvexity, SCP and DBHM converge to different solutions even applied with the same initial point.
Because the subproblems are solve separately, different subproblems consider different wind conditions. In this way, the optimal layout is explored
in a larger scope, which makes DBHM more possible to find a better result than direct SCP solution.

\begin{figure}[h]
	\centering
		\includegraphics[width=3.2 in]{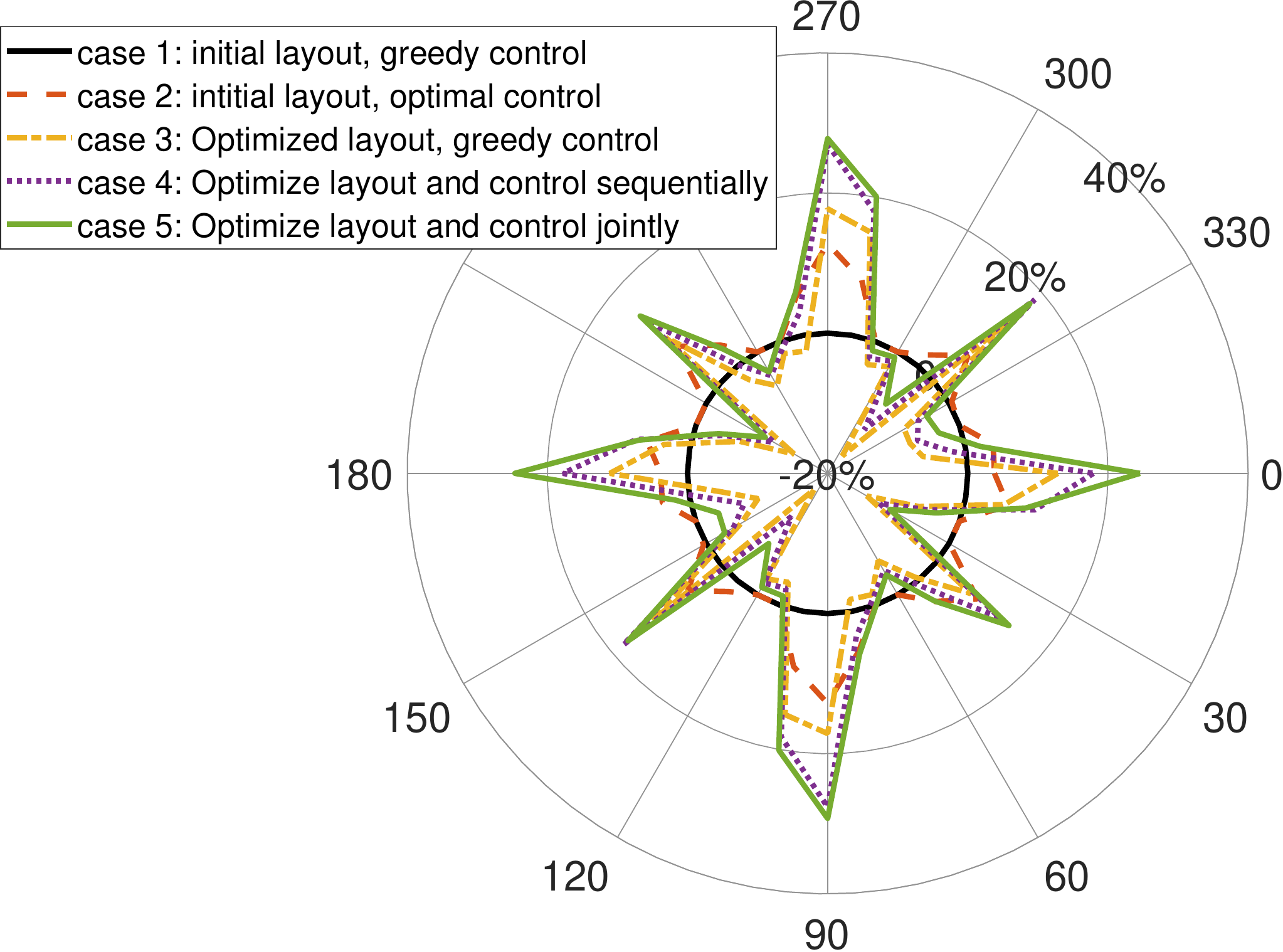}
  \vspace{-5pt}
	\caption{AEP increase with respect to case 1 for different wind directions.}
	\label{fig:EPPimprovement}
\end{figure}

Fig. \ref{fig:EPPimprovement} depicts the power improvement with respect to case 1 for different wind directions. For both the initial and optimized layouts, the control optimization can always increase the power production. Compared with the initial rectangular layout, the gains of layout optimization vary with the wind direction. The increase can reach up to $20\%$-$30\%$ for the wind directions parallel to the rows of the initial layout, such as $0^{o}$, $90^{o}$, $180^{o}$, and $270^{o}$. For directions such as $30^{o}$, $130^{o}$, $210^{o}$ and $310^{o}$, the power production is sacrificed, but an average increase is obtained for the overall wind distribution. 

To make the results more stable, we also consider finer wind direction scenarios. Table \ref{Tab360} shows the results with 360 wind direction scenarios which is fine enough to evaluate the AEP \cite{feng2015solving}. Compared to the optimization results of 36 directions, the overall trends of different cases are similar. The joint optimization still obtains the best result. The AEP is increased by $1.44\%$ compared to the sequent optimization. On the other hand, the AEP improvement by the control optimization increases. This is because the control optimization considers individual wind condition separately and the finer direction discretization increases the optimization space. On the contrary, evaluated in 360 scenarios, the AEP improvements by the layout optimization are smaller than those evaluated in fewer scenarios, which aligns well with the layout optimization analysis in the previous research \cite{feng2015solving}. This can be explained by the fact that in the layout optimization, turbine placements are manipulated considering the average power production for the overall wind directions. The power production in different wind directions need to be traded-off. The low direction resolution simplifies the problem by assuming that wind only comes from fewer specific directions, which leaves a large space for optimization algorithms to adjust the WF layout. With finer wind direction discretization, the layout optimization needs to balance the WF power production in more different scenarios and the optimization margin is thus smaller.

\begin{table}[h]
\centering
\caption{AEP Results of Various Optimization Cases for the 16-WT WF with 360 Scenarios}
\label{Tab360}
\begin{tabular}{c c c c c}
\hline
  Case & Layout & Operation & AEP(GWh) & Improvement \\
\hline
  1 & Initial & Greedy &   374.34 &   -- \\
  2 & Initial & Optimized & 385.40 &  2.95\% \\
  3 & Optimized & Greedy & 378.34 &   1.07\% \\
  4 & Optimized & Optimized (sequent) &  386.39 &  3.22\% \\
  5 & Optimized & Optimized (Joint) &  391.97 &  4.71\% \\
\hline
\end{tabular}
\end{table}

\subsection{WF with larger turbine distance}
\begin{figure}[h]
	\centering
		\includegraphics[width=2.9 in]{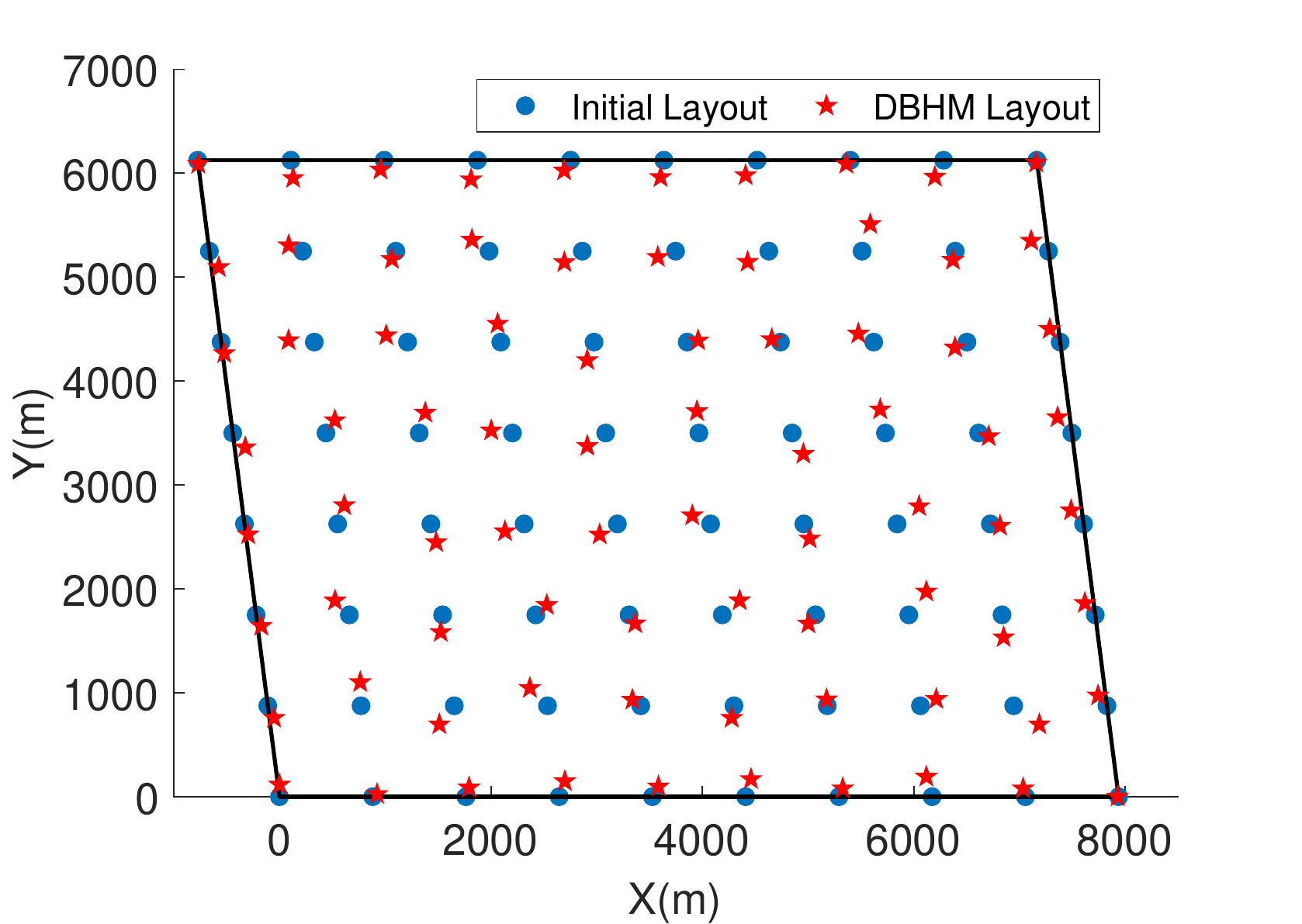}
	\caption{Initial and optimized layout of the 80-WT WF.}
	\label{fig:80layout}
\end{figure}

In this subsection, we proceed to a sparser WF where the distances between the turbines are much larger. In such sparser WFs, the wake interactions are weaker. 
The initial layout is shown in Fig. \ref{fig:80layout}. Eighty WTs replicates the Horns Rev1 WF in Denmark \cite{park2015layout}. Turbines are placed in the corners of a parallelogram layout with a $7.2^{o}$ tilt, whose side lengths are approximately 7D. The wind direction distribution is shown in Fig \ref{fig:windrose}. The investigated cases are analogue to the cases in subsection \ref{16wtcase}. The AEP improvements of different optimization cases with respect to case 1 are compared in Table \ref{Tab4}.  Compared to the case in subsection \ref{16wtcase}, the Horn Rev1 WF is sparser, and the wake interactions are weaker. Therefore, the optimization improvement is lower. In this case, wind direction is discretized into 12 scenarios. The total number of decision variables is 2080. For the decomposition formulation, 12 subproblems are obtained with 320 variables each. 
Even in this sparser WF, DBHM can still improve the AEP by $0.95 \%$ compared to the sequential optimization of case 4 by exploiting the synergy between WF layout and operation.

\begin{table}[h]
\centering
\caption{AEP Results of Various Optimization Cases for the 80-WT WF}
\label{Tab4}
\begin{tabular}{c c c c c}
\hline
  Case & Layout & Operation & AEP(GWh) & Improvement \\
\hline
  1 & Initial & Greedy &   1971.1 &   -- \\
  2 & Initial & Optimized &   1984.3 &  0.67\% \\
  3 & Optimized & Greedy &   2041.0 &   3.54\% \\
  4 & Optimized & Optimized (sequent) &   2058.4 &   4.40\% \\
  5 & Optimized & Optimized (Joint) &  2078.0 &  5.42\% \\
\hline
\end{tabular}
\end{table}

For the isolated WF layout optimization by PSO in case 3, the computational time is 300629 s. Compared to the 16-WT case, the computational burden increase rapidly with the increase of variable dimensions as analyzed in subsection \ref{comptime}. For the joint optimization (\ref{obj})-(\ref{wakeconstraint}), the computational burden of heuristic methods is clearly infeasible. Table \ref{Tab5} compares the computation time and results when solving the joint optimization model with different methods. Leveraging the decomposability of the joint optimization problem, DBHM obtains the best result for this considered joint optimization problem.

\begin{table}[h]
\centering
\caption{Computational Performance of Different Methods}
\label{Tab5}
\begin{tabular}{c c c c}
\hline
  Method & SCP & PSO & DBHM \\
\hline
  Average Computation Time(s) & 280220 & $>$ 10 days &  314188  \\
  AEP(GWh) & 2058.93 & - &   2078.0  \\
  Improvement compared to case 4 & 0.02\% & - &  0.95 \% \\
\hline
\end{tabular}
\end{table}

%

The finer wind direction resolution with 180 scenarios is also considered in Table \ref{Tabmore80WT}. Compared to the $12$-scenario case, different AEP results are obtained and the differences among cases 2-5 decrease. Evaluated in 180 scenarios, the AEP of the initial case 1 is higher. Besides, as analyzed in the 16-WT WF, the improvements by the layout optimization are smaller when evaluated in 180 scenarios, because the WF power production in more different wind directions are considered and balanced. The AEP improvement by the control optimization increases. On the other hand, the overall trends of different optimization cases do not change. Both the layout and the control optimization can improve the AEP. Besides, the joint optimization still obtains the best result. Compared to the sequent optimization, the AEP is improved by $0.38\%$. The optimization results using 180 scenarios is much more reliable. In Fig. \ref{fig:80layout}, the optimized WF layout obtained by the joint optimization considering 180 scenarios is also shown. 

\begin{table}[h]
\centering
\caption{AEP Results of Various Optimization Cases for the 80-WT WF with 180 Scenarios}
\label{Tabmore80WT}
\begin{tabular}{c c c c c}
\hline
  Case & Layout & Operation & AEP(GWh) & Improvement \\
\hline
  1 & Initial & Greedy &   1981.1 &   -- \\
  2 & Initial & Optimized &   1998.1 &  0.86\% \\
  3 & Optimized & Greedy &  1999.9 &   0.95\% \\
  4 & Optimized & Optimized (sequent) & 2015.2 &   1.72\% \\
  5 & Optimized & Optimized (Joint) & 2022.9 &  2.11\% \\
\hline
\end{tabular}
\end{table}

\section{Conclusion}\label{sec:cons}
In this paper, the WF layout optimization is studied with optimal cooperative control considerations. A novel two-stage WF layout optimization model is proposed so that the control optimization can be considered in the layout design phase to obtain the system-level optimal layout design. The WF layout is co-optimized with the WF operation and their synergy are exploited. Besides, a DBHM method is proposed, which can leverage the hierarchy and decomposability of the joint optimization problem and solve the problem efficiently.

By comparing the WF layouts obtained by the proposed joint optimization with the layouts obtained by traditional separate optimizations, the synergy between WF layout and operation is analyzed. Two WFs are studied. For a 16-WT WF, the AEP is increased by $1.44\%$ by doing the joint optimization compared to the separate optimization. Besides, it is found that for the cooperative WF control, the AEP increase mainly comes from the non-zero yaw angles. Therefore, for the joint optimization, the axial induction factors can be fixed at 1/3. Only the yaw and the layout are required to be co-optimized for AEP maximization. For a 80-WT WF with larger turbine distance, by considering the cooperative control in the layout optimization phase, the AEP can still be improved by $0.38\%$. In all cases, the joint optimization obtains the best result, which shows that the synergy between the layout and the operation is beneficial to be considered in the WF layout design. 

For the novel joint optimization problem considered in this paper, DBHM obtains better results with high computational efficiency compared with applying traditional layout optimization methods directly. 

In this work, we focus on the optimization of power extraction. The proposed method will be further improved to consider more complex objectives in the future, such as the transportation cost and the maintenance cost.



%

\section*{Appendix A} \label{DBHM}
\subsection*{Overall Iteration Framework of the Proposed DBHM}
\begin{figure}[h]
	\centering
		\includegraphics[width=3.2 in]{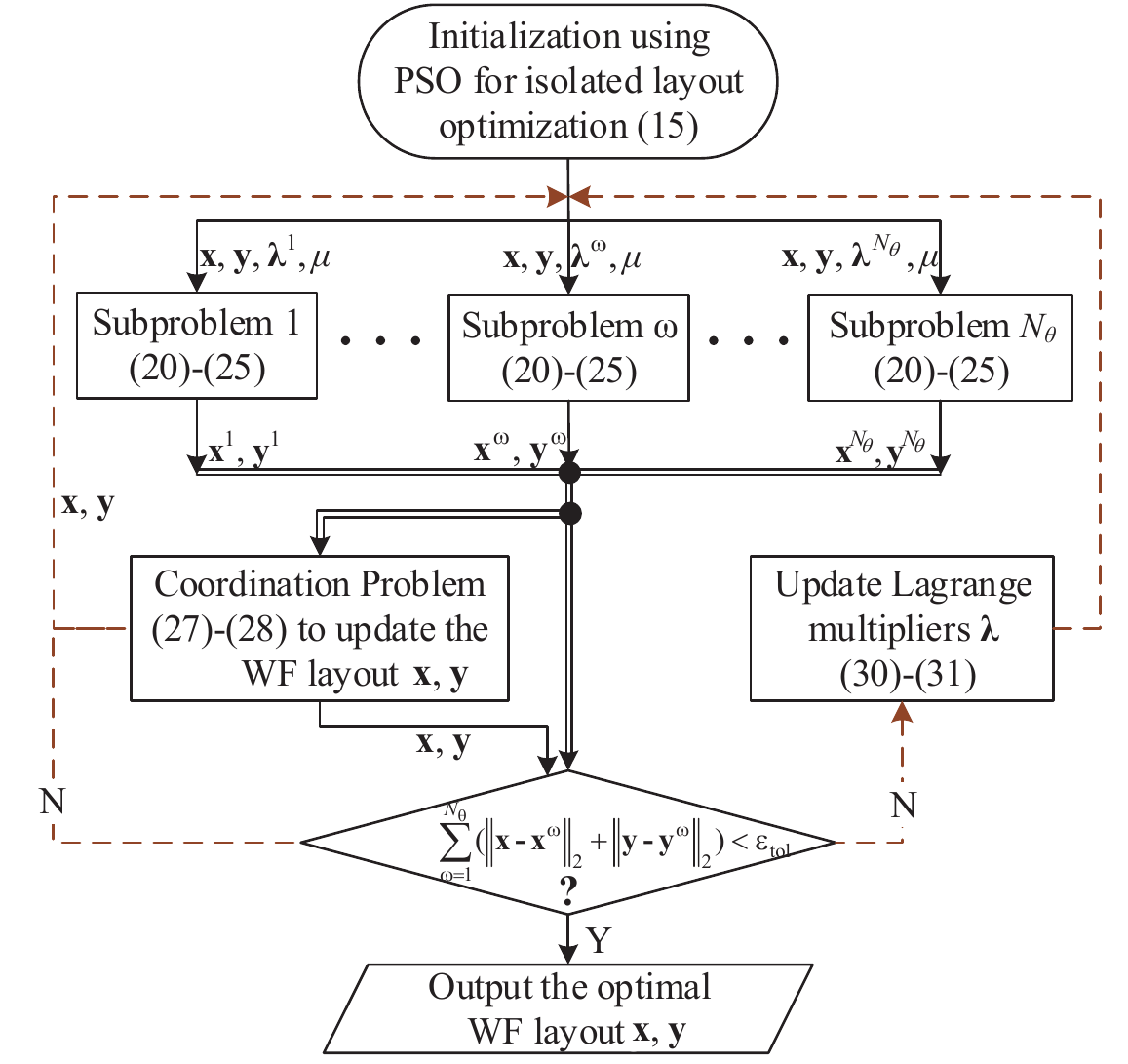}
	\caption{Flowchart of the proposed DBHM.}  
	\label{fig:Flowchart}
\end{figure}

The overall flowchart of the proposed decomposition-based hybrid method (DBHM) for the joint WF layout and control optimization is shown in Fig. \ref{fig:Flowchart}.
As shown in Fig. \ref{fig:Flowchart}, the heuristic PSO is first used to solve the isolated layout optimization to provide a good initial point. Then, the subproblems (\ref{subobj})-(\ref{multip}) and their coordination problem (\ref{upperx})-(\ref{uppery}) are solved iteratively following the process of alternating direction method of multipliers \cite{boyd2011distributed, wang2019global}. At each iteration, the individual subproblems of all wind scenarios are first solved either sequentially or in parallel given $\textbf{x},\textbf{y}$ and $\lambda, \mu$, each of which considers only one particular wind condition. Then, under the same $\lambda$ and $\mu$, (\ref{upperx}) and (\ref{uppery}) are computed to update the WF layout variables $\textbf{x}$ and $\textbf{y}$ by coordinating different wind scenarios. Because the subproblem is much smaller than the original joint optimization model, the computational burden is reduced and can be distributed. 

The iteration terminates when the violation of consensus constraints (\ref{admmconsesus}) is lower than the pre-defined tolerance $\epsilon_{\textrm{tol}}$
\begin{align}
  \sum_{\omega=1}^{N_\theta} (\| \textbf{x} - \textbf{x}^{\omega}  \|_2 + \| \textbf{y} - \textbf{y}^{\omega}  \|_2) < \epsilon_{\textrm{tol}}
\end{align}

If the convergence criterion is not satisfied, the Lagrange multipliers $\lambda$ are updated and the iteration is performed again.

Under the augmented Lagrangian formulation of the subproblem (\ref{subobj})-(\ref{multip}), $\lambda$ is updated based on the violation of the corresponding consensus constraints \cite{boyd2011distributed}:
\begin{align}
   \lambda_{i,x}^{\omega,k+1} = \lambda_{i,x}^{\omega,k} +  2*\mu (x_i - x_i^{\omega}), \hspace{8pt} \forall \omega, i \\
   \lambda_{i,y}^{\omega,k+1} = \lambda_{i,y}^{\omega,k} +  2*\mu (y_i - y_i^{\omega}), \hspace{8pt} \forall \omega, i
\end{align}
where the superscript $k$ indicates the values at the $k$th iteration.
The step size is the corresponding penalty factor $\mu$. This updating tries to force $\lambda_{i,x}$ and $\lambda_{i,y}$ to converge to the Lagrange multipliers of a local optimum of the original problem (\ref{admobj})-(\ref{admmope}), so that the WF layout $\textbf{x}$ and $\textbf{y}$ can consequently converge to the local optimum of the original problem \cite{boyd2011distributed}.

For sufficiently large $\mu$, the iteration process shown in Fig. \ref{fig:Flowchart} is guaranteed to converge to the set of stationary solutions of the original joint optimization model (\ref{obj})-(\ref{wakeconstraint}). The convergency of the consensus optimization with compact feasible sets has been theoretically proved by the recent study of decomposition theories  \cite{wang2019global}. The proposed decomposition formulation (\ref{admobj})-(\ref{admmope}) fulfills the two requirements and thus the convergency can be guaranteed. Firstly, the feasible region of the joint optimization model restricted by constraints (\ref{poscons})-(\ref{terrcons2}) and (\ref{opecons})-(\ref{opecons2}) is a compact set, because the layout area and control actions are always bounded in realistic WF problems. Secondly, the coupling constraints between the introduced local variables ($\textbf{x}^{\omega},\textbf{y}^{\omega}$) and original variables ($\textbf{x},\textbf{y}$) are the linear consensus constraints (\ref{admmconsesus}). For this form of nonconvex optimizations, according to \cite{wang2019global}, the iteration process of alternating direction method of multipliers converges to stationary solutions of the original optimization model as long as the penalty factor $\mu$ is chosen large enough. In this paper, $\mu = 10$ which is adjusted and can guarantee the convergency of the proposed model. We omit the tedious convergence theories of the nonconvex decomposition problem, which is not the contribution of this paper. The detailed mathematical proof can be found in \cite{wang2019global}. 

\section*{References}

  \bibliographystyle{elsarticle-num}
  \bibliography{bare_jrnl}

\begin{thebibliography}{10}
\expandafter\ifx\csname url\endcsname\relax
  \def\url#1{\texttt{#1}}\fi
\expandafter\ifx\csname urlprefix\endcsname\relax\def\urlprefix{URL }\fi
\expandafter\ifx\csname href\endcsname\relax
  \def\href#1#2{#2} \def\path#1{#1}\fi

\bibitem{councilwind}
G.~W. E.~C. GWEC, Global offshore wind report 2020, URL
  https://gwec.net/global-offshore-wind-report-2020/.

\bibitem{veers2019grand}
P.~Veers, K.~Dykes, Lantz, et~al., Grand challenges in the science of wind
  energy, Science 366~(6464) (2019) eaau2027.

\bibitem{boersma2017tutorial}
S.~Boersma, B.~M. Doekemeijer, P.~M. Gebraad, P.~A. Fleming, J.~Annoni, A.~K.
  Scholbrock, J.~A. Frederik, J.-W. van Wingerden, A tutorial on
  control-oriented modeling and control of wind farms, in: 2017 American
  control conference (ACC), IEEE, 2017, pp. 1--18.

\bibitem{ahmadyar2016coordinated}
A.~S. Ahmadyar, G.~Verbi{\v{c}}, Coordinated operation strategy of wind farms
  for frequency control by exploring wake interaction, IEEE Trans on
  Sustainable Energy 8~(1) (2016) 230--238.

\bibitem{dar2016windfarm}
Z.~Dar, K.~Kar, O.~Sahni, J.~H. Chow, Windfarm power optimization using yaw
  angle control, IEEE Transactions on Sustainable Energy 8~(1) (2016) 104--116.

\bibitem{van2019effects}
D.~van~der Hoek, S.~Kanev, J.~Allin, D.~Bieniek, N.~Mittelmeier, Effects of
  axial induction control on wind farm energy production-a field test,
  Renewable Energy 140 (2019) 994--1003.

\bibitem{fleming2017field}
P.~Fleming, J.~Annoni, J.~J. Shah, L.~Wang, S.~Ananthan, Z.~Zhang,
  K.~Hutchings, P.~Wang, W.~Chen, L.~Chen, Field test of wake steering at an
  offshore wind farm, Wind Energy Science Discussions 2~(NREL/JA-5000-67623).

\bibitem{park2016bayesian}
J.~Park, K.~H. Law, Bayesian ascent: A data-driven optimization scheme for
  real-time control with application to wind farm power maximization, IEEE
  Transactions on Control Systems Technology 24~(5) (2016) 1655--1668.

\bibitem{hou2019review}
P.~Hou, J.~Zhu, K.~Ma, G.~Yang, W.~Hu, Z.~Chen, A review of offshore wind farm
  layout optimization and electrical system design methods, Journal of Modern
  Power Systems and Clean Energy 7~(5) (2019) 975--986.

\bibitem{GAO2016192}
X.~Gao, H.~Yang, L.~Lu, Optimization of wind turbine layout position in a wind
  farm using a newly-developed two-dimensional wake model, Applied Energy 174
  (2016) 192--200.

\bibitem{tao2020wind}
S.~Tao, Q.~Xu, A.~Feij{\'o}o, G.~Zheng, J.~Zhou, Wind farm layout optimization
  with a three-dimensional gaussian wake model, Renewable Energy.

\bibitem{tao2019optimal}
S.~Tao, S.~Kuenzel, Q.~Xu, Z.~Chen, Optimal micro-siting of wind turbines in an
  offshore wind farm using frandsen--gaussian wake model, IEEE Trans on Power
  Systems 34~(6) (2019) 4944--4954.

\bibitem{gonzalez2013new}
J.~S. Gonz{\'a}lez, M.~B. Pay{\'a}n, J.~R. Santos, A new and efficient method
  for optimal design of large offshore wind power plants, IEEE Transactions on
  Power Systems 28~(3) (2013) 3075--3084.

\bibitem{bastankhah2016experimental}
M.~Bastankhah, F.~Port{\'e}-Agel, Experimental and theoretical study of wind
  turbine wakes in yawed conditions, Journal of Fluid Mechanics 806 (2016) 506.

\bibitem{annoni2018analysis}
J.~Annoni, P.~Fleming, A.~K. Scholbrock, J.~M. Roadman, S.~Dana, C.~Adcock,
  F.~Porte-Agel, S.~Raach, F.~Haizmann, D.~Schlipf, Analysis of
  control-oriented wake modeling tools using lidar field results, Wind Energy
  Science (Online) 3~(NREL/JA-5000-72767).

\bibitem{wu2021design}
C.~Wu, X.~Yang, Y.~Zhu, On the design of potential turbine positions for
  physics-informed optimization of wind farm layout, Renewable Energy 164
  (2021) 1108--1120.

\bibitem{long2017formulation}
H.~Long, Z.~Zhang, Z.~Song, A.~Kusiak, Formulation and analysis of grid and
  coordinate models for planning wind farm layouts, IEEE Access 5 (2017)
  1810--1819.

\bibitem{cao2020optimizing}
J.~F. Cao, W.~J. Zhu, W.~Z. Shen, J.~N. S{\o}rensen, Z.~Y. Sun, Optimizing wind
  energy conversion efficiency with respect to noise: A study on multi-criteria
  wind farm layout design, Renewable Energy 159 (2020) 468--485.

\bibitem{turner2014new}
S.~Turner, D.~Romero, P.~Zhang, C.~Amon, T.~Chan, A new mathematical
  programming approach to optimize wind farm layouts, Renewable Energy 63
  (2014) 674--680.

\bibitem{MITTAL2016133}
P.~Mittal, K.~Kulkarni, K.~Mitra, A novel hybrid optimization methodology to
  optimize the total number and placement of wind turbines, Renewable Energy 86
  (2016) 133--147.

\bibitem{reddy2020wind}
S.~R. Reddy, Wind farm layout optimization (windflo): An advanced framework for
  fast wind farm analysis and optimization, Applied Energy 269 (2020) 115090.

\bibitem{brogna2020new}
R.~Brogna, J.~Feng, J.~N. S{\o}rensen, W.~Z. Shen, F.~Port{\'e}-Agel, A new
  wake model and comparison of eight algorithms for layout optimization of wind
  farms in complex terrain, Applied Energy 259 (2020) 114189.

\bibitem{PEREZ2013389}
B.~Pérez, R.~Mínguez, R.~Guanche, Offshore wind farm layout optimization
  using mathematical programming techniques, Renewable Energy 53 (2013)
  389--399.

\bibitem{park2015layout}
J.~Park, K.~H. Law, Layout optimization for maximizing wind farm power
  production using sequential convex programming, Applied energy 151 (2015)
  320--334.

\bibitem{feng2015solving}
J.~Feng, W.~Z. Shen, Solving the wind farm layout optimization problem using
  random search algorithm, Renewable Energy 78 (2015) 182--192.

\bibitem{2012Wind}
C.~Wan, J.~Wang, G.~Yang, H.~Gu, X.~Zhang, Wind farm micro-siting by gaussian
  particle swarm optimization with local search strategy, Renewable Energy 48
  (2012) 276--286.

\bibitem{ju2019wind}
X.~Ju, F.~Liu, Wind farm layout optimization using self-informed genetic
  algorithm with information guided exploitation, Applied Energy 248 (2019)
  429--445.

\bibitem{huang20183}
L.~Huang, H.~Tang, K.~Zhang, Y.~Fu, Y.~Liu, 3-d layout optimization of wind
  turbines considering fatigue distribution, IEEE Transactions on Sustainable
  Energy 11~(1) (2018) 126--135.

\bibitem{YAMANIDOUZISORKHABI2018341}
S.~{Yamani Douzi Sorkhabi}, D.~A. Romero, J.~C. Beck, C.~H. Amon, Constrained
  multi-objective wind farm layout optimization: Novel constraint handling
  approach based on constraint programming, Renewable Energy 126 (2018)
  341--353.

\bibitem{wilson2018evolutionary}
D.~Wilson, S.~Rodrigues, C.~Segura, I.~Loshchilov, F.~Hutter, G.~L. Buenfil,
  A.~Kheiri, E.~Keedwell, M.~Ocampo-Pineda, E.~{\"O}zcan, et~al., Evolutionary
  computation for wind farm layout optimization, Renewable energy 126 (2018)
  681--691.

\bibitem{EROGLU201253}
Y.~Eroğlu, S.~U. Seçkiner, Design of wind farm layout using ant colony
  algorithm, Renewable Energy 44 (2012) 53--62.

\bibitem{long2015two}
H.~Long, Z.~Zhang, A two-echelon wind farm layout planning model, IEEE
  Transactions on Sustainable Energy 6~(3) (2015) 863--871.

\bibitem{hou2016offshore}
P.~Hou, W.~Hu, M.~Soltani, C.~Chen, B.~Zhang, Z.~Chen, Offshore wind farm
  layout design considering optimized power dispatch strategy, IEEE Trans on
  sustainable energy 8~(2) (2016) 638--647.

\bibitem{zong2021experimental}
H.~Zong, F.~Port{\'e}-Agel, Experimental investigation and analytical modelling
  of active yaw control for wind farm power optimization, Renewable Energy 170
  (2021) 1228--1244.

\bibitem{wang2016novel}
L.~Wang, A.~Tan, Y.~Gu, A novel control strategy approach to optimally design a
  wind farm layout, Renewable Energy 95 (2016) 10--21.

\bibitem{pedersen2020integrated}
M.~M. Pedersen, G.~C. Larsen, Integrated wind farm layout and control
  optimization, Wind Energy Science 5~(4) (2020) 1551--1566.

\bibitem{fleming2016wind}
P.~A. Fleming, A.~Ning, P.~M. Gebraad, K.~Dykes, Wind plant system engineering
  through optimization of layout and yaw control, Wind Energy 19~(2) (2016)
  329--344.

\bibitem{doekemeijer2019tutorial}
B.~M. Doekemeijer, J.-W. Van~Wingerden, P.~A. Fleming, A tutorial on the
  synthesis and validation of a closed-loop wind farm controller using a
  steady-state surrogate model, in: 2019 American Control Conference (ACC),
  IEEE, 2019, pp. 2825--2836.

\bibitem{2016Maximization}
P.~Gebraad, J.~J. Thomas, A.~Ning, P.~Fleming, K.~Dykes, Maximization of the
  annual energy production of wind power plants by optimization of layout and
  yaw based wake control, Wind Energy 20~(1) (2016) 97--107.

\bibitem{fathy2001coupling}
H.~K. Fathy, J.~A. Reyer, P.~Y. Papalambros, A.~Ulsov, On the coupling between
  the plant and controller optimization problems, in: Proceedings of the 2001
  American Control Conference.(Cat. No. 01CH37148), Vol.~3, IEEE, 2001, pp.
  1864--1869.

\bibitem{zhong2016decentralized}
S.~Zhong, X.~Wang, Decentralized model-free wind farm control via discrete
  adaptive filtering methods, IEEE Transactions on Smart Grid 9~(4) (2016)
  2529--2540.

\bibitem{reddy2021efficient}
S.~R. Reddy, An efficient method for modeling terrain and complex terrain
  boundaries in constrained wind farm layout optimization, Renewable Energy 165
  (2021) 162--173.

\bibitem{2010Runtime}
D.~Sudholt, C.~Witt, Runtime analysis of a binary particle swarm optimizer,
  Theoretical Computer ence 411 (2010) 2084--2100.

\bibitem{hou2015optimized}
P.~Hou, W.~Hu, M.~Soltani, Z.~Chen, Optimized placement of wind turbines in
  large-scale offshore wind farm using particle swarm optimization algorithm,
  IEEE Trans on Sustainable Energy 6~(4) (2015) 1272--1282.

\bibitem{POOKPUNT2013266}
S.~Pookpunt, W.~Ongsakul, Optimal placement of wind turbines within wind farm
  using binary particle swarm optimization with time-varying acceleration
  coefficients, Renewable Energy 55 (2013) 266--276.

\bibitem{boyd2011distributed}
S.~Boyd, N.~Parikh, E.~Chu, Distributed optimization and statistical learning
  via the alternating direction method of multipliers, Now Publishers Inc,
  2011.

\bibitem{jonkman2009definition}
J.~Jonkman, S.~Butterfield, W.~Musial, G.~Scott, Definition of a 5-mw reference
  wind turbine for offshore system development, National Renewable Energy
  Laboratory, Golden, CO, Technical Report No. NREL/TP-500-38060.

\bibitem{campagnolo2020wind}
F.~Campagnolo, R.~Weber, J.~Schreiber, C.~L. Bottasso, Wind tunnel testing of
  wake steering with dynamic wind direction changes, Wind Energy Science 5~(4)
  (2020) 1273--1295.

\bibitem{E2017Statistical}
E.~Thøgersen, B.~Tranberg, J.~Herp, M.~Greiner, Statistical meandering wake
  model and its application to yaw-angle optimisation of wind farms, Journal of
  Physics: Conference Series 854 (2017) 012017--.

\bibitem{munters2018dynamic}
W.~Munters, J.~Meyers, Dynamic strategies for yaw and induction control of wind
  farms based on large-eddy simulation and optimization, Energies 11~(1) (2018)
  177.

\bibitem{2016Wind}
F.~Campagnolo, V.~Petrovic, C.~L. Bottasso, A.~Croce, Wind tunnel testing of
  wake control strategies, in: American Control Conference, 2016.

\bibitem{2010Global}
C.~Mathworks, Global optimization toolbox user ' s guide.

\bibitem{wang2019global}
Y.~Wang, W.~Yin, J.~Zeng, Global convergence of admm in nonconvex nonsmooth
  optimization, Journal of Scientific Computing 78~(1) (2019) 29--63.

\end{thebibliography}

\nolinenumbers
\end{document}